\pdfoutput=1

\documentclass[a4paper,11pt]{amsart}

\usepackage{graphicx}
\usepackage{amsmath,amscd,amsfonts,amssymb,latexsym}
\usepackage[cmtip,arrow]{xy}
\usepackage{pb-diagram,pb-xy}
\usepackage{graphicx}
\usepackage{tikz}
\usepackage{url}

\newcommand{\C}{\mathbb{C}}

\newcommand{\T}{\mathbb{T}}

\newcommand{\Z}{\mathbb{Z}}
\newcommand{\PP}{\mathbb{P}}
\newcommand{\G}{G}
\newcommand{\Gdwa}{{\bf G}_2}
\newcommand{\HH}{{\rm H}}

\newcommand{\uu}{\mathcal U}
\newcommand{\Hom}{{\rm Hom}}
\newcommand{\dlog}{{\rm dlog}}
\newcommand{\oA}{{A^o}}

\newtheorem{theorem}{Theorem}
\newtheorem{proposition}[theorem]{Proposition}

\newtheorem{remark}[theorem]{Remark}

\title{Residues formulas for the push-forward in K-theory, the case of $\Gdwa/P$.}
\author{ Andrzej Weber \and Magdalena Zielenkiewicz}
\thanks{The first author supported by NCN grant 2013/08/A/ST1/00804. }\thanks{The second author supported by NCN grant 2015/17/N/ST1/02327.}
\address{University of Warsaw, Institute of Mathematics \\
 Banacha 2, 02-097 Warszawa, Poland}
 \email{aweber@mimuw.edu.pl}
\address{University of Warsaw, Institute of Mathematics\\
 Banacha 2, 02-097 Warszawa, Poland}
 \email{magdaz@mimuw.edu.pl}

\begin{document}

\maketitle

\begin{abstract}
We study residue formulas for push-forward in the equivariant K-theory of homogeneous spaces. For the classical Grassmannian the residue formula can be obtained from the cohomological formula by a substitution. We also give another proof using symplectic reduction and the method involving the localization theorem of Jeffrey--Kirwan. We review formulas for other classical groups, and we derive them from the formulas for the classical Grassmannian. Next we consider the homogeneous spaces for the exceptional group $\Gdwa$. One of them, $\Gdwa/P_2$ corresponding to the shorter root, embeds in the Grassmannian $Gr(2,7)$. We find its fundamental class in the equivariant K-theory $K^\T(Gr(2,7))$. This allows to derive a residue formula for the push-forward. It has significantly different character comparing to the classical groups case. The factor involving the fundamental class of $\Gdwa/P_2$ depends on the equivariant variables. To perform computations more efficiently we apply the basis of K-theory consisting of Grothendieck polynomials. The residue formula for push-forward remains valid for the homogeneous space $\Gdwa/B$ as well.
\end{abstract}

\section{Introduction}
\label{intro}
Suppose an algebraic torus $\T=(\C^*)^r$ acts on a complex variety $X$ which is smooth and complete. Let $E$ be an equivariant complex vector bundle over $X$. Then $E$ defines an element of the equivariant K-theory of $X$. In our setting there is no difference whether one considers the algebraic K-theory or the topological one. We study the push-forward $p_!$ of the class of $E$ to the equivariant K-theory of a point (for the basics of equivariant K-theory see e.g. \cite{Chriss})
$$K^\T(pt)\simeq {\rm R}(\T)\simeq\Z[t_1,t_1^{-1},t_2,t_2^{-1},\dots,t_r,t_r^{-1}]\,.$$
Suppose that the fixed point set $X^\T$ is finite. The localization formula allows to compute the push-forward using calculus of rational functions:
\begin{equation}\label{locformula}p_!(E)=\sum_{x\in X^\T}\frac{E_x}{\prod_{i=1}^{d} \left(1-\tfrac1{L_i(x)}\right)}\,.\end{equation}
Here the $T_xX=\bigoplus_{i=1}^{d}L_i(x)$ is the decomposition of the tangent space into a sum of one-dimensional representations. Identifying the representation ring ${\rm R}(\T)$ with Laurent polynomials, one obtains a sum whose summands are rational functions in $t_i$. After simplification the result is a Laurent polynomial in $t_i$.

Formula (\ref{locformula}) is a K-theoretic version of the cohomological Atiyah--Bott--Berline--Vergne localization formula \cite{AtBo2,BeVe}. It has appeared in the literature much earlier, see e.g.
\cite{AtBo1}, \cite[\S4.7]{Nielsen},
\cite[Cor.~6.12]{Groth}, \cite{BFQ}.
This formula is also called ``Holomorphic Lefschetz Theorem'' in \cite[(4.6)]{AtiyahSinger}. A relative version of the localization formula is given in \cite[Th. 5.11.7]{Chriss}.
In some of the quoted papers one considers the actions of finite groups but the modification to a torus action is purely formal.
\medskip

We study the case when $X$ is a homogeneous space with the maximal torus action. One of the examples is the Grassmannian $X=Gr(m,n)$ of the $m$-dimensional spaces in $\C^n$. Suppose that $E$ is a vector bundle associated to a representation of $GL_m(\C)$. Let $f(z_1,\dots,z_m)$ be the character of this representation. Then the formula (\ref{locformula}) takes the form
\begin{equation}\label{locformulagrass}p_!(E)=\sum_{A\subset\{t_1,t_2,\dots,t_n\},\,|A|=m}\; \frac{f(A)}{
\prod_{a\in A}\prod_{b\in \oA}\left(1-\frac{a}{b}\right)}\,.\end{equation}
where
$\oA=  \{t_1,t_2,\dots,t_n\}\setminus A$ and $f(A)=f(t_{i_1},t_{i_2},\dots,t_{i_m})$
for a subset $A=\{t_{i_1},t_{i_2},\dots,t_{i_m}\}$ of the basic characters. (In section \S\ref{notacja} we introduce notation which allows to write this formula in a compact manner.)
\medskip

Suppose $X$ is a space obtained by the symplectic reduction construction. This is the case for the Grassmannian, since
$$Gr(m,n)=\Hom(\C^m,\C^n)/\!\!/U(m)\,.$$ Equivalently, $Gr(m,n)$ is a GIT quotient
$\Hom(\C^m,\C^n)/\!\!/GL_m(\C)\,.$ Precisely, the Grassmannian is the geometric quotient of the open set consisting of the injective maps $\C^m\to \C^n$. Then the push-forward in equivariant K-theory can be expressed by a residue, just as in the case of cohomology. A generalization of the theory of Jeffrey--Kirwan \cite{JeKi} or Guillemin--Kalkman \cite{GuKa} applied to the homogeneous spaces leads to explicit residue formulas for the push-forward in equivariant cohomology. Such formulas for complete flag varieties first appeared in \cite{BeSz}. They were generalized to Grassmannians for the classical groups in \cite{Zie1}. A conceptual proof based on the Jeffrey--Kirwan theory is presented in \cite{Zie2}.
Similar formulas can be obtained in equivariant K-theory.
Formulas for Grassmannians and complete flag varieties appeared already in \cite{All}, see also \cite{AllRi}.
In \cite{All} the proof is omitted since it is stated that the proof is similar to \cite{Zie1} and that it just ``requires
more notation and paper''. Nevertheless  the proofs can be explained conceptually.
The Jeffrey--Kirwan method can be applied for a class of generalized cohomology theories as shown in \cite{Metz}.
For the Grassmannian $Gr(m,n)$ the formula for the push-forward obtained via symplectic reduction is the following
\begin{equation}\label{residueforgrass}p_!(E)=\frac1{m!}Res_{Z=0,\infty}\frac{f(Z)\prod_{z\not=z'}(1-\tfrac{z}{z'})}{\prod_{t\in T}\prod_{z\in Z}(1-\tfrac zt)}\prod_{z\in Z}\frac{dz}z\,,\end{equation}
where $Z=\{z_1,z_2,\dots,z_m\}$ is an auxiliary set of variables.
The operation $Res_{Z=0,\infty}$ is the iterated residue at $0$ and $\infty$
$$(Res_{z_1=0}+Res_{z_1=\infty})\circ(Res_{z_2=0}+Res_{z_2=\infty})\circ\dots\circ(Res_{z_m=0}+Res_{z_m=\infty})\,.$$
In our case the order of taking the residues does not matter.
\medskip

The formula for the classical Grassmannian allows one to write immediately the residue formulas for orthogonal and Lagrangian Grassmannians, using the embeddings $OG(2n)\subset Gr(n,2n)$, $OG(2n+1)\subset Gr(n,2n+1)$ and $LG(n)\subset Gr(n,2n)$. The results are presented in section \S\ref{odgrassmanianu}. The general case requires embedding the homogeneous space in question into a flag variety and presenting the flag variety as a symplectic reduction.  For classical groups this is achieved in \cite{Zie3}, where the cohomology formulas are derived. The case of K-theory is analogous and is discussed shortly in \S\ref{drugidowod}, and we give examples of resulting formulas in the final section. Our main goal is to present residue formulas for the exceptional group $\Gdwa$. An embedding of $\Gdwa/P_2$ (see \S\ref{G2opis} for the notation and \cite{Anthesis} for the description of homogeneous spaces of $\Gdwa$) into $Gr(2,7)=\Hom(\C^2,\C^7)/\!\!/GL_2(\C)$ is known. But the inverse image of $\Gdwa/P_2$ in $\Hom(\C^2,\C^7)$ is not a complete intersection, as is the case of classical groups. Moreover, its class in equivariant K-theory (also in equivariant cohomology) is significantly more complicated.
The residue formula computing push-forward contains a factor which is given by Theorem \ref{g2alpha}. The case of the other homogeneous space $\Gdwa/P_1$ is simpler, since this space embeds as a quadric in $\PP^6$.
\bigskip

We would like to thank the anonymous Referees for many important comments which helped to improve the paper essentially.

\section{Notation}\label{notacja}
Let $\T$ be a torus. Let $\T^\vee=\Hom(\T,\C^*)$. The equivariant K-theory of a point is isomorphic to the group algebra
$$K^\T(pt) \simeq \Z[\T^\vee].$$
When $\T\simeq (\C^*)^n$, then the group algebra $K^\T(pt)$ can be identified with Laurent polynomials  $\Z[t_1^{\pm1},t_2^{\pm1},\dots,t_n^{\pm1}]$. The generator $t_i$ denotes the character (the representation) $\T\to \C^*=GL_1(\C)$, which is the projection on the $i$-th factor. We use the multiplicative notation for characters. For example the character $t_2t_4$ denotes the character $$\T \to\C^*\,,\qquad
(t_1,t_2,\dots,t_n)\mapsto t_2t_4\,,$$
while $t_2+t_4$ stands for two dimensional representation which decomposes into two one dimensional representations having characters $t_2$ and $t_4$.
The unit 1 is the class of the trivial representation of dimension 1.

\subsection{Operations on sets of characters}

For a set of characters $\{a_1,a_2,\dots,a_k\}=A\subset \T^\vee$ let

\begin{itemize}
\item $A^{-1}=\{a^{-1}\;|\;a\in A\}$

\item $\Lambda(A)=\{a_ia_j\;|\;1\leq i<j\leq k\}$

\item $S(A)=\{a_ia_j\;|\;1\leq i\leq j\leq k\}$

\item $R(A)=\{\tfrac ab\;|\;a,b\in A,\;a\not=b\}$

\item $R^+(A)=\{a_i/a_j\;|\;i<j\}.$
\end{itemize}
For two sets of characters $A,B\subset \T^\vee$ let
\begin{itemize}
\item $A/B=\{\tfrac ab\;|\;a\in A,\;b\in B\}$
%\item $A\otimes B=\{ab\;|\;a\in A,\;b\in B\}$
\end{itemize}
The definition of $R^+(A)$ depends on the ordering of $A$ while the remaining ones do not.
Note that if $A$ is the set of coordinate characters $t_i$ for the maximal torus in $GL_n(\C)$, then  $R(A)$ is the set of roots and $R^+(A)$ is the set of positive roots.

\subsection{Meromorphic functions associated to sets of characters}

For a set of characters $A\subset \T^\vee$ let
\begin{itemize}
\item $[A]=\prod_{a\in A}(1-\frac1a)\in K^\T(pt)$
\item $\dlog(A)=\prod_{a\in A}\frac {da}a$
\end{itemize}
This expressions are understood as meromorphic functions or  forms on the product of $\PP^1$'s.

\subsection{Notation for multi-variable residue}

For a meromorphic form $f(z)dz$ on $\PP^1$ let $$Res_{z=0,\infty}f\,dz=Res_{z=0}f\,dz+Res_{z=\infty}f\,dz\,.$$
For a meromorphic form $f(z_1,z_2,\dots,z_n)dz_1dz_2\dots dz_n$, possibly depending on additional parameters, the following operator
$$Res_{z_1=0,\infty}\circ Res_{z_{2}=0,\infty}\circ\dots\circ Res_{z_n=0,\infty}$$
will be denoted by $Res_{Z=0,\infty}$.

\section{The standard Grassmannian}\label{gra1}
Denote by $Gr(m,n)$ the Grassmannian of $m$-subspaces in $\C^n$.
Let $E$ be a vector bundle over the Grassmannian $Gr(m,n)$ which is associated to a representation of $GL_m(\C)$. Let $f\in K^\T(pt)$ be the associated character.
It is a symmetric function $f(A)=f(t_1,t_2,\dots ,t_m)$  in the coordinate characters of $\T=(\C^*)^m$.
Let
\begin{itemize}
\item $T=\{t_1,t_2,\dots, t_n\}$,
\item $Z=\{z_1,z_2,\dots ,z_m\}$,
\item $\oA=T\setminus A$ for $A\subset T$.
\end{itemize}
Let $p:{Gr(m,n)}\to pt$ and let $p_!:K^\T({Gr(m,n)})\to K^\T(pt)$ be the push-forward in the equivariant K-theory. In the introduction we have described a residue formula for the push-forward in K-theory of the Grassmannian.
In the notation introduced in \S\ref{notacja} the push-forward for the classical Grassmannian can be expressed as
\begin{equation}\label{standard}p_!(E)=\sum_{A\subset T,\;|A|=m}\frac{ f(A)}{[\oA/A]}=
\frac{1}{m!} Res_{Z=0,\infty} \frac{f(Z)\cdot[R(Z)]\cdot \dlog(Z)}{[T/Z]}\,.\end{equation}
The expression above is nothing but the formulas (\ref{locformulagrass}) and (\ref{residueforgrass}) written in our notation.
Another (more compact) formula is given by:
$$p_!(E)= Res_{Z=0,\infty} \frac{f(Z)\cdot[R^+(Z)]\cdot \dlog(Z)}{[T/Z]}.$$
%=Res_{Z=0,\infty} \frac{f(Z)\cdot \dlog(Z)}{[T/Z]}\,.$$
Here there is no factor $\frac1{m!}$, and the numerators are simpler: one can take only positive roots $R^+(Z)$, without altering the result. A proof of this formula is sketched in \cite[Prop 5.3]{All}.

\section{Proof of residue formulas}
The  formula (\ref{standard}) can be deduced from the cohomological push-forward formula (see \cite{We}, Chapter 4. or \cite{Zie1}, Corollary 3.2) or can be proved directly using symplectic reduction as in \cite{Zie2,Zie3}.
\subsection{Reduction to the cohomological formula}\label{dowody}

For a polynomial $g$ in $Z$ one has
$$\sum_{A\subset T,\;|A|=m}\frac{ g(A)}{\prod_{x\in A\, t\in\oA}(t-x)}
=\frac1{m!}Res_{Z=0}\frac{ g(Z)\prod_{z\not=z'}(z-z') dZ}{\prod_{t\in T\, z\in Z}(t-z)}.$$
The same formula is valid for Laurent polynomials, provided that we take into account the residues at zero.
Let $g(Z)= f(Z)/\prod_{z\in Z}z^{m}$. Then
$$\sum_{A\subset T,\;|A|=m}\frac{ f(A)}{\prod_{t\in A}t^{m}\prod_{x\in A\, t\in\oA}(t-x)}
=\frac1{m!}Res_{Z=0,\infty}\frac{ f(Z)\prod_{z\not=z'}(z-z') dZ}{\prod_{z\in Z}z^{m}\prod_{t\in T\, z\in Z}(t-z)}$$
The left hand side is equal to
$$\sum_{A\subset T,\;|A|=m}\frac{ f(A)}{\prod_{t\in T}t^{m}\prod_{x\in A\, t\in\oA}(1-x/t)}\,,$$
while the right hand side equals
\begin{multline*}\frac1{m!}Res_{Z=0,\infty}\frac{ f(Z)\prod_{z\not=z'}(1-z'/z)\prod_{z\in Z}z^{m-1} dZ}{\prod_{z\in Z}z^{m}\prod_{t\in T\, z\in Z}(t-z)}=\\=\frac1{m!}Res_{Z=0,\infty}\frac{ f(Z)\prod_{z\not=z'}(1-z'/z)\prod_{z\in Z}z^{-1} dZ}{\prod_{t\in T}t^m\prod_{t\in T\, z\in Z}(1-z/t)}\end{multline*}
Multiplying both sides by $\prod_{t\in T}t^{m}$
 we obtain the formula (\ref{standard}), that is
$$\sum_{A\subset T,\;|A|=m}\frac{ f(A)}{\prod_{x\in A\, t\in\oA}(1-x/t)}=\frac1{m!}Res_{Z=0,\infty}\frac{ f(Z)\prod_{z\not=z'}(1-z'/z)\prod_{z\in Z}z^{-1} dZ}{\prod_{t\in T\, z\in Z}(1-z/t)}\,.$$
\qed

\medskip

Note that if $f$ is homogeneous of degree at least $m$ in each variable, then the residue at $0$ vanishes.

\subsection{Proof via symplectic reduction}\label{drugidowod}

The proof  is based on the fact that $Gr(m,n)$ is a symplectic reduction of $V=\Hom(\C^m,\C^n)$ and the following two propositions:

\begin{proposition} \label{prop:jk} For a representation $$[E]\in 
K^{GL_m(\C)\times \T}(V)\simeq K^{GL_m(\C)\times \T}(pt)$$ there is an equality
$$p_!( E_{V/\!\!/GL_m(\C)})=q_!( \varpi E_{V/\!\!/\T_m})\in K^\T(pt)\,.$$
Here $\T_m$ is a maximal torus in $GL_m(\C)$, $q:V/\!\!/\T_m\to pt$ and
$$\varpi=\prod_{\alpha\in \text{ roots of } GL_m(\C)}(1-\alpha).$$\end{proposition}
The factor $\varpi$ is equal to $[R(Z)]$ according to our notation, provided that $z_i$'s are the standard coordinates of the diagonal maximal torus in $GL_m(\C)$. The symplectic reduction $V/\!\!/\T_m$ is a product of projective spaces.

\begin{proposition} \label{prop:proj} For a class $E\in K^{\C^*\times \T}(\C^n)\simeq K^{\C^*\times \T}(pt)$ given by the character $$f(z,t_1,t_2,\dots, t_n)$$
we have
$$q_!(E_{\PP^{n-1}})=Res_{z=0,\infty}\left(\frac{f(z,t_1,t_2,\dots, t_n)}{z\prod_{i=1}^n(1-\tfrac z{t_i})}dz\right)$$
\end{proposition}

Proposition \ref{prop:jk} is a K-theoretic version of the Jeffrey--Kirwan Theorem (\cite{JeKi}). For the proof see \cite{harada}, Theorem 5.3.

\medskip
\noindent{\it Proof of Proposition \ref{prop:proj}.} The map
$$f(z,t_1,t_2,\dots, t_n)\mapsto Res_{z=0,\infty}\left(\frac{f(z,t_1,t_2,\dots, t_n)}{z\prod_{i=1}^n(1-\tfrac z{t_i})}dz\right)$$
is linear with respect to the ring $\Z[t_1^{\pm1},t_2^{\pm1},\dots,t_n^{\pm1}]$, hence it is enough to examine the formula of the proposition for $f(z,t_1,t_2,\dots, t_n)=z^k$. Let us first consider the case  $k=-\ell\leq 0$. Then the form
$$g(z)=\frac{f(z,,t_1,t_2,\dots, t_n)}{z\prod_{i=1}^n(1-\tfrac z{t_i})}dz=\frac{dz}{z^{\ell+1}\prod_{i=1}^n(1-\tfrac z{t_i})}$$ has no pole at infinity and the residue at 0 is equal to the coefficient of $z^\ell$ in the expansion
$$\frac{1}{\prod_{i=1}^n(1-\tfrac z{t_i})}=\prod_{i=1}^n\sum_{j=0}^\infty \tfrac {z^j}{t_i^j}\,.$$
This coefficient is equal to the $\ell$--th complete symmetric function in the inverses of variables $t_i$
$$\sum_{|\lambda|=\ell}t^{-\lambda }\,.$$
This is the character of the representation $Sym^\ell(\C^n)^*=q_!(\mathcal O(\ell))$, that is exactly what we had to show.

For $f=z^k$, $k\in[1,n-1]$ the function $g(z)$ has no pole at 0 nor at infinity. For $k\geq n$ we only have a residue at infinity. Substituting $z$ by $w^{-1}$ and applying the transformation
$$\frac{z^k}{z\prod_{i=1}^n(1-\tfrac z{t_i})}dz=
\frac{-1}{w^{(k-1)}\prod_{i=1}^n(1-\tfrac 1{wt_i})}\frac{dw}{w^2}=
\frac{(-1)^{n-1}\prod_{i=1}^nt_i}{w^{k+1-n}\prod_{i=1}^n(1-wt_i)}dw$$
we obtain the result by  Serre duality.
\section{The residue formulas derived from \S\ref{gra1}}
\label{odgrassmanianu}
Below we consider Lagrangian and orthogonal Grassmannians which are invariant with respect to smaller tori. In all the considered cases the fixed point set in the big Grassmannian with respect to the small torus is equal to the fixed point set with respect to the full torus action. Therefore we can apply the residue formula and specialize it to the smaller torus. In all the cases we apply the projection formula: for $\iota:G\hookrightarrow Gr(k,n)$ we have
$$(p_{|G})_!(\iota^*E)=p_!([G]\cdot E)\,,$$
where $[G]=\iota_!(1)$ is the fundamental class of $G$ in the K-theory of the Grassmannian.
\subsection{Lagrangian Grassmannian}
Let $LG(n)\subset Gr(n,2n)$ be the Lagrangian Grassmannian. Set
\begin{itemize}
\item $T=\{t_1,t_2,\dots, t_n\}$,
\item $T^{\pm1}=T\cup T^{-1}$,
\item $Z=\{z_1,z_2,\dots ,z_n\}$.
\end{itemize}
Let  $p_L:LG(n)\to pt$. Then

$$p_{L_!}(E)=\sum_{A\subset T}\frac{ f(A\cup \oA{^{-1}})}{[S(A^{-1}\cup \oA)]}=
\frac{1}{n!} Res_{Z=0,\infty} \frac{f(Z)\cdot[\Lambda (Z^{-1})]\cdot[R(Z)]\cdot \dlog(Z)}{[T^{\pm1}/Z]}\,.$$
The first equality follows is the Atiyah--Bott--Berline--Vergne formula for $LG(n)$. The fixed points are indexed by subsets of $T$ and the tangent space at the fixed subspace $W_A$ corresponding to the subset $A$ is canonically identified with $Sym^2(W_A^{*})$. The second equality is obtained by applying the formula for the classical Grassmannian to the embedding of $LG(n)$ in $Gr(n, {2n})$ and restricting to the action of a smaller torus via the embedding $(t_1,\dots, t_n, t_n^{-1}, \dots, t_1^{-1}) \mapsto (t_1,\dots, t_{2n})$. The factor $\Lambda(Z)$ corresponds to the fundamental class of the embedding, as follows. The set of Lagrangian maps $\C^n\to \C^{2n}$ is a complete intersection in $\Hom(\C^n,\C^{2n})$, described by the $\frac{n(n-1)}{2}$ independent equations of the form $\omega(v_i, v_j)=0$ for $1 \leq i < j \leq n$. Here $\omega$ is the symplectic form and $v_i$ are the column vectors of the matrix corresponding to an element of $\Hom(\C^n,\C^{2n})$. The set of weights of the equations is equal to $\Lambda(Z^{-1})$, where $z_i$ are characters associated to the source.

\subsection{Orthogonal Grassmannian, even case}
Let $OG(2n)\subset Gr(n,{2n})$ be the isotropic Grassmannian (two components). As before, set
\begin{itemize}
\item $T=\{t_1,t_2,\dots, t_n\}$,
\item $T^{\pm1}=T\cup T^{-1}$,
\item $Z=\{z_1,z_2,\dots ,z_n\}$.
\end{itemize}
Let  $p_O:OG(2n)\to pt$. Then

$$p_{O!}(E)=\sum_{A\subset T}\frac{ f(A\cup \oA{^{-1}})}{[\Lambda(A^{-1}\cup \oA)]}=
\frac{1}{n!} Res_{Z=0,\infty} \frac{f(Z)\cdot[S( Z^{-1})]\cdot[R(Z)]\cdot \dlog(Z)}{[T^{\pm1}/Z]}\,.$$
As in the case of the Lagrangian Grassmannian, the first equality is simply the Atiyah--Bott--Berline--Vergne formula. The fixed points are indexed by subsets of $T$, and the tangent space at the subspace $W_A$ corresponding to a subset $A$ is canonically identified with $\Lambda^2(W^*_A)$. Note that the set of isotropic maps $\C^n\to \C^{2n}$ is a complete intersection in $\Hom(\C^n,\C^{2n})$, described by the $\frac{n(n+1)}{2}$ independent equations of the form $\Omega(v_i, v_j)=0$ for $1 \leq i \leq j \leq n$. Here $\Omega$ is a nondegenerate symmetric form and $v_i$ are the column vectors of the matrix corresponding to an element of $\Hom(\C^n,\C^{2n})$. The set of weights of the equations is equal to $S(Z^{-1})$, where $z_i$ are characters associated to the source.

\subsection{Orthogonal Grassmannian, odd case}

Let $OG(2n + 1)\subset Gr(n,{2n+1})$ be the isotropic Grassmannian. As before, set
\begin{itemize}
\item $T=\{t_1,t_2,\dots, t_n\}$,
\item $T^{\sharp}=T\cup T^{-1}\cup\{1\}$,
\item $Z=\{z_1,z_2,\dots ,z_n\}$.
\end{itemize}
Let  $p_O:OG(2n+1)\to pt$. Then

\begin{multline*}p_{O!}(E)=\sum_{A\subset T,\;|A|}\frac{ f(A\cup \oA{^{-1}})}{[\Lambda(A^{-1}\cup \oA)]\cdot [A^{-1}\cup \oA]} =\\
\frac{1}{n!} Res_{Z=0,\infty} \frac{f(Z)\cdot[S( Z^{-1})]\cdot[R(Z)] \cdot \dlog(Z)}{[T^{\sharp}/Z]}\,.
\end{multline*}
The set of isotropic maps $\C^n\to \C^{2n+1}$ is a complete intersection in $\Hom(\C^n,\C^{2n})$, described by the same equations as for the even-dimensional case. The set of weights of the equations is therefore equal to $S(Z^{-1})$, where $z_i$ are characters associated to the source. Note that here we have a factor $[z_i^{-2}]=1-z_i^2$ in the numerator and $[z_i^{-1}]=1-z_i$ in the denominator. We simplify the whole expression and obtain
$$p_{O!}(E)
=
\frac{1}{n!} Res_{Z=0,\infty} \frac{f(Z)\cdot[\Lambda( Z^{-1})]\cdot[-Z^{-1}]\cdot[R(Z)] \cdot \dlog(Z)}{[T^{\pm1}/Z]}\,.
$$

\section{Grassmannian of $\phi$-isotropic planes}
Consider the root system of the exceptional group $\Gdwa$, with simple roots $\alpha_1$ and $\alpha_2$.
 
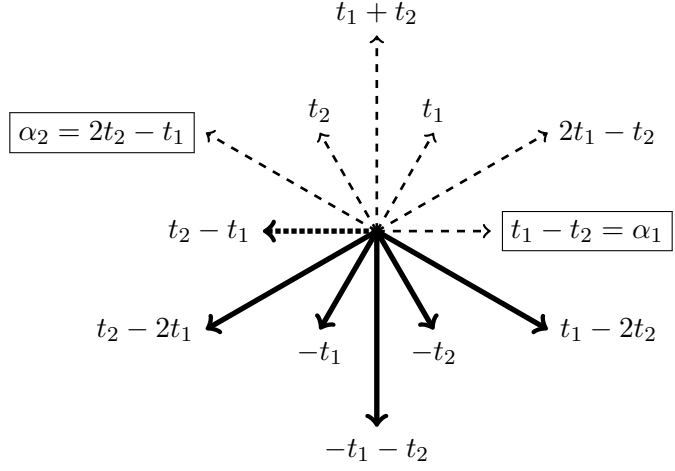
\begin{center}
\begin{tikzpicture}[scale=1.5]
	\draw[line width=1,color=black,dashed,->] (0,0) -- (1,0) node[anchor=west]{$\boxed{ t_1-t_2=\alpha_1}$};
%        \path (1,0) node[anchor=west]
	\draw[line width=2,color=black,densely dotted,->] (0,0) -- (-1,0) node[anchor=east]{$t_2-t_1$};
	\draw[line width=1,color=black,dashed,->] (0,0) -- (0,1.73) node[anchor=south]{$ t_1+t_2$};
	\draw[line width=2,color=black,->] (0,0) -- (0,-1.73) node[anchor=north]{$-t_1-t_2$};
	\draw[line width=1,color=black,dashed,->] (0,0) -- (3/2,0.87) node[anchor=west]{$ 2t_1-t_2$};
	\draw[line width=2,color=black,->] (0,0) -- (3/2,-0.87) node[anchor=west]{$t_1-2t_2$};
	\draw[line width=1,color=black,dashed,->] (0,0) -- (-3/2,0.87) node[anchor=east]{$\boxed{\alpha_2={2t_2-t_1}}$};
	\draw[line width=2,color=black,->] (0,0) -- (-3/2,-0.87) node[anchor=east]{$t_2-2t_1$};
	\draw[line width=1,color=black,dashed,->] (0,0) -- (1/2,0.87) node[anchor=south]{$t_1$};
	\draw[line width=2,color=black,->] (0,0) -- (1/2,-0.87) node[anchor=north]{$-t_2$};
	\draw[line width=1,color=black,dashed,->] (0,0) -- (-1/2,0.87) node[anchor=south]{$ t_2$};
	\draw[line width=2,color=black,->] (0,0) -- (-1/2,-0.87) node[anchor=north]{$-t_1$};
\end{tikzpicture}
\end{center}
\begin{center}\footnotesize The positive roots -- dashed arrows,\\ the only negative root of $P_2 \supset B$ -- the bold dotted root $\alpha_1$,\\ the characters of the tangent space $T_{[1]}\Gdwa/P_2$ -- bold arrows.\end{center}

\noindent Let $B$ be the Borel subgroup defined by the positive roots.
The group $\Gdwa$ has exactly two maximal parabolic subgroups containing $B$, corresponding to the two choices of simple roots. We denote by $P_i$ the parabolic subgroup in which the  root $-\alpha_i$ is omitted.
 Anderson in \cite{An} gave an interpretation of homogeneous spaces for $\Gdwa$. Given a nondegenerate trilinear form $\phi: \Lambda^3 \C^7 \to \C$, a two dimensional subspace $V$ is called $\phi$-isotropic if $\phi(x,y,z)=0$ for all $x,y \in V$, $z\in\C^7$.
The homogeneous space $\Gdwa/P_2$ parameterizes two-dimensional $\phi$-isotropic subspaces in $\C^7$, yielding an embedding of $\Gdwa/P_2$ in the Grassmannian $Gr(2,7)$.
The embedding allows us to derive a residue formula for push-forward.

\begin{theorem}\label{g2alpha}
Suppose a bundle $E\to \Gdwa/P_2$ is associated to a character $f(Z)=f(z_1,z_2)$. Then
$$p_{\Gdwa!}(E)= Res_{Z=0,\infty}\frac{f(Z)\cdot\uu(Z,T)\cdot [R^+(Z)]\cdot\dlog(Z)}{[T^\flat/Z]}\,,$$
where $$T^\flat=\{t_1 , t_2 , t_1t_2^{-1}, 1, t_1^{-1} t_2 , t_2^{-1} , t_1^{-1}\}$$
and
$$ \uu(Z,T)= {z_1 z_2}\left(1 - {z_1}\right) \left(1 - {z_2}\right) \left(1 - {z_1 z_2}\right) \left(\uu_z(Z)-\uu_t(T)\right)\,,$$
with
$$\uu_z(Z)=z_1+\frac{1}{z_1}+z_2+\frac{1}{z_2}+z_1 z_2+\frac{1}{z_1 z_2}-6\,,$$
$$\uu_t(T)=\uu_z\left(t_1,\frac1{t_2}\right)=t_1+\frac{1}{t_1}+t_2+\frac{1}{t_2}+\frac{t_1}{t_2}+\frac{t_2}{t_1}-6\,.$$
\end{theorem}  
We describe the embedding $\Gdwa/P_2\hookrightarrow Gr(2,7)$ in detail and derive the residue formula in the next section. For the details of computations see the Appendix.
Also there it is given 
a formula obtained directly from Atiyah-Bott-Berline-Vergne theorem. That is a sum of six rational functions, see the formula (\ref{cyclicsum}).
\section{Computation of the polynomial $\uu(Z,T)=[\Gdwa/P_2]\subset K^T(Gr(2,7))$}
\label{G2opis}

 The exceptional Lie algebra $\mathfrak g_2$ is the fixed set of the triality automorphism of $\mathfrak{so}(8)$. The triality automorphism of the real form  $\mathfrak{so}(4,4)$ permutes the off-diagonal entries (up to a sign). The resulting embedding
 $\mathfrak{g}_2\hookrightarrow\mathfrak{gl}_8$ computed in \cite{MiWe} is the following:
\begin{equation}\label{zanurzenie}\left(
\begin{array}{cc|cccccc}
 {t_1} & {a_1} & -{a_2} & {a_3} &
   {a_3} & {a_4} & {a_5} & 0 \\
 {a_6} & {t_2} & {a_7} & {a_2} & {a_2} &
   {a_3} & 0 & -{a_5} \\
\hline
\bf -{a_8} &\bf {a_9} & {t_1-t_2} & {a_1} & {a_1}
   & 0 & -{a_3} & -{a_4} \\
 \bf{a_{10}} & \bf{a_8} & {a_6} & 0 & 0 & -{a_1} &
   -{a_2} & -{a_3} \\
 \bf{a_{10}} & \bf{a_8} & {a_6} & 0 & 0 & -{a_1} &
   -{a_2} & -{a_3} \\
 \bf{a_{11}} & \bf{a_{10}} & 0 & -{a_6} & -{a_6} &
   {-t_1+t_2} & -{a_7} & {a_2} \\
 \bf{a_{12}} & 0 & -{a_{10}} & -{a_8} & -{a_8} &
   -{a_9} & -{t_2} & -{a_1} \\
 0 & \bf-{a_{12}} & -{a_{11}} & -{a_{10}} & -{a_{10}} &
  {a_8} & -{a_6} & -{t_1}
\end{array}
\right)\end{equation}
In this presentation of $\mathfrak{g}_2\subset\mathfrak{gl}_8$  the 4-th and the 5-th coordinates are equal, hence $\mathfrak{g}_2$ annihilates the vector $(0,0,0,1,-1,0,0,0)$ and therefore preserves its orthogonal complement. Hence
$$\mathfrak{g}_2\subset\mathfrak{gl}({\rm span}\{\varepsilon_1,\varepsilon_2,\varepsilon_3,\varepsilon_4+\varepsilon_5,\varepsilon_6,\varepsilon_7,\varepsilon_8\})\simeq\mathfrak{gl}_7\,.$$
The parabolic group $P_2$ is equal to the intersection of $\Gdwa$ with the standard block-upper-triangular matrices in $GL_7(\C)$ preserving the first two coordinates.
The weights of the tangent space $T_{[1]}\Gdwa/P_2$ are the following: \newline

\begin{center}
\begin{tabular}{r|c|c|c|c|c}
coordinate &
$a_8$&
$a_9$&
$a_{10}$&
$a_{11}$&
$a_{12}$\\ \hline
weight&
 $-t_2$& $t_1-2t_2$ & $-t_1$ & $t_2-2t_1$& $-t_1-t_2$
\end{tabular}
\end{center}

\medskip

Another description of the embedding is given in \cite{Anthesis}.
Given a two-dimensional space $W$ one  identifies  $\C^7$ with $W \oplus End^0(W) \oplus W^{*}$, where $End^0(W)$ denotes the trace $0$ endomorphisms of $W$.  The antisymmetric trilinear form $\phi$ on $\C^7$ is induced by the map
$$W \times End^0(W) \times V^{*}\to \C\,,$$
$$(v,F,f)\mapsto f(F(v))$$
and $W$ is isotopic with respect to this form.
The tangent space to $[W]$ in $\Gdwa/P_2$ coincides, via this embedding, with
$$ (Sym^3 W^{*} \otimes \Lambda^2 W) \oplus \Lambda^2 W^{*}\hookrightarrow  \Hom(W, End^0(W) \oplus W^{*}).$$
The weights of the action of the tangent space (in the multiplicative notation) are therefore $\{ t_1^{-2}t_2\,,\, t_1^{-1}\,,\, t_2^{-1}\,,\, t_1 t_2^{-2}\,,\, t_1^{-1}t_2^{-1} \}$. This agrees with the embedding (\ref{zanurzenie}).

Since $\Gdwa/P_2$ embeds in the Grassmannian $Gr(2,7)$, the push-forwards in K-theory are related as follows. Let $\T$ denote the maximal torus in $\Gdwa$. The set of characters obtained from the embedding to $GL_7(\C)$ described by (\ref{zanurzenie}) is equal to $$T^\flat=\{t_1 \,,\, t_2 \,,\, t_1t_2^{-1}\,,\, 1\,,\, t_1^{-1} t_2 \,,\, t_2^{-1} \,,\, t_1^{-1}\}\,.$$ The torus $\T$ acts on $Gr(2,7)$ leaving $\Gdwa/P_2$ invariant.
Given a class $a\in K^\T(Gr(2,7))$, let $\iota^* a$ denote its restriction to $K^\T({\Gdwa/P_2})$ via the embedding $$\iota:\Gdwa/P_2 \hookrightarrow Gr(2,7)\,.$$
Then $p_{\Gdwa!}(\iota^*a)= p_!([\Gdwa/P_2] \cdot  a)$, where $$p_!: K^\T(Gr(2,7)) \to K^\T(pt) \quad\text{and}\quad  p_{\Gdwa!}: K^\T({\Gdwa/P_2}) \to K^\T(pt)$$ are the respective push-forwards. If the class
 $\iota^* a$ is represented by a bundle $E\to \Gdwa/P_2$ associated to a character $f(Z)=f(z_1,z_2)$, then by the residue formula (\ref{standard}) for $Gr(2,7)$ one gets
$$p_{\Gdwa!}(E)= Res_{Z=0,\infty}\frac{f(Z)\cdot\uu(Z,T)\cdot [R^+(Z)]\cdot\dlog(Z)}{[T^\flat/Z]}\,,$$
where
$$\uu(Z,T)\in  K^{ GL(2)\times\T}(pt)\subset \Z[z_1,z_1^{-1},z_2,z_2^{-1},t_1,t_1^{-1},t_2,t_2^{-1}]$$ is a lift of the fundamental class $[\Gdwa/P_2]\subset K^\T(Gr(2,7))$.
The lift $\uu(Z,T)$ is not unique and might be quite complicated.
It remains to find a convenient form of $\uu(Z,T)$. We have found a formula which seems to have a nice and remarkable shape.
\medskip

We start with a warm-up  presenting the computation of the cohomology class of $\Gdwa/P_2$ in $\HH^*(Gr(2,7))$. It is simpler than finding the class in K-theory, but essentially the method is the same. The computations presented in this section have been performed using Wolfram Mathematica, see the Appendix and \cite{Wolfram}.
To find the cohomology class  we integrate the Schur classes $S_{5}$, $S_{41}$ and $S_{32}$ of the dual of the tautological bundle. These classes form a self-dual basis of  $\HH^{10}(Gr(2,7))$. We find (e.g. using Atiyah--Bott--Berline--Vergne integration formula) that
$$\int_{\Gdwa/P_2}S_{5}(\gamma^*)=0\,,\qquad \int_{\Gdwa/P_2}S_{41}(\gamma^*)=2\,,\qquad\int_{\Gdwa/P_2}S_{32}(\gamma^*)=2\,.$$
Therefore
\begin{equation}\label{homologicznie}[\Gdwa/P_2]=2S_{41}(\gamma^*)+2S_{32}(\gamma^*)\in \HH^{10}(Gr(2,7))\,.\end{equation}
\medskip

An analogous calculation can be performed in equivariant cohomology. The equivariant Schur classes $S_I=S_I(\gamma^*)$, for $I$ contained in the rectangle $5\times 2$, form a basis of $\HH^*_\T(Gr(2,7))$ as a free module over $\HH^*_\T(pt)$. 
We note that the Schur polynomials we use here are not the double Schur polynomials (representing equivariant classes of Schubert varieties) but we apply usual Schur polynomials to the equivariant bundle $\gamma^*$.
Then, except from the classes  $S_{41}$ and $S_{32}$, we obtain non-zero result for
$$\int_{\Gdwa/P_2}S_{52}(\gamma^*)=4 (t_1^2 - t_1 t_2 + t_2^2)\,,\qquad
\int_{\Gdwa/P_2}S_{43}(\gamma^*) = 2 (t_1^2 - t_1 t_2 + t_2^2)\,,$$
$$\int_{\Gdwa/P_2}S_{54}(\gamma^*) = 2 (t_1^2 - t_1 t_2 + t_2^2)^2\,.$$
Inverting the intersection matrix we obtain a coefficients in the equivariant Schur basis:
$$[\Gdwa/P_2]=2S_{41}+2S_{32}- 2(t_1^2 - t_1 t_2 + t_2^2)S_{21}\in \HH^{10}_\T(Gr(2,7))\,.$$
This shows that the equivariant fundamental class of $\Gdwa/P_2$ has higher order coefficients in the Schur basis. Such situation does not appear for the classical groups.
In the terms of the Chern roots $x_1,x_2$ of the bundle $\gamma^*$ the fundamental class is equal to
$$2 x_1 x_2 (x_1 + x_2) [(x_1^2 + x_1 x_2 + x_2^2)
- (t_1^2 - t_1 t_2 + t_2^2)]\,.$$
\medskip

The same procedure leads to finding the class in K-theory. There is a small difference, which in fact is not a problem from the theoretical point of view: there is no grading in K-theory. One has to take into account all partitions contained in the $5\times 2$ rectangle. We replace the Schur basis by the basis consisting of the Grothendieck classes of $\gamma^*$, i.e.~the (stable) Grothendieck polynomials in the formal roots of $\gamma^*$ (\cite{LaSch,Buch}). Grothendieck classes represent the  classes of the structure sheaves of Schubert varieties in non-equivariant K-theory. Our notation agrees with \cite{Buch,RS}. This basis has ``skew-upper-triangular'' intersection form, which means that
$$p_!(\G_I(\gamma^*)\G_J(\gamma^*))=0\in K(pt)=\Z\quad\text{if}\quad I\not\subset c(J)\,.$$
Here $c(J)$ denotes the complementary partition. The intersections of Schubert varieties are rational varieties, therefore one has
$$p_!(\G_I(\gamma^*)\G_J(\gamma^*))=1\quad\text{if}\quad I\subset c(J)\,.$$
We compute the push-forward of the (nonequivariant) Grothendieck classes $\G_I(\gamma^*)$ for the map $\Gdwa/P_2\to pt$, e.g. applying the localization formula \cite[Cor.~6.12]{Groth}. We obtain the following list of results:
$$
\begin{array}{c|c|c|c|c|c|c|c|c|c|c|c|c}
\text{partition } I &[0] & [1] & [2] & [11] & [3] &
   [21] & [40] & [31] & [22] & [50] &
   [41] & [32]  \\ \hline
 p_{\Gdwa !}(S_I(\gamma^*))~
&1 & 1 & 1 & 1 & 1 & 1 & 2 & 1 & 1 & 0 & 2 & 2  \\
\end{array}
$$
The push-forwards for the partitions $I$ with $|I|>5$ vanish since for those partitions $\G_I$ is supported by a variety of the codimension $>5=\dim \Gdwa/P_2$.
Inverting the intersection matrix of $K(Gr(2,7))$ we find that the nonequivariant class of $\Gdwa/P_2$ in the K-theory is equal to
\begin{equation}\label{G2nonequiv}2 \G_{41}+2 \G_{32} - \G_{33}  - 3 \G_{4 2} + \G_{4 3}\,.\end{equation}
(We omit $\gamma^*$ in the notation.)
This class can be written in the Grothendieck roots of the tautological bundle as
$$  {z_1 z_2}\left(1 - {z_1}\right) \left(1 - {z_2}\right) \left(1 - {z_1 z_2}\right)
\left(z_1+\frac{1}{z_1}+z_2+\frac{1}{z_2}+z_1 z_2+\frac{1}{z_1 z_2}-6\right)
\,.$$
%$$  \frac1{z_1 z_2}\left(1 - \frac1{z_1}\right) \left(1 - \frac1{z_2}\right) \left(1 - \frac1{z_1 z_2}\right)
%\left(z_1+\frac{1}{z_1}+z_2+\frac{1}{z_2}+z_1 z_2+\frac{1}{z_1 z_2}-6\right)
%\,.$$
In the same way one computes the equivariant class. The equivariant Grothendieck classes of the bundle $\gamma^*$ form a basis of $K^T(\Gdwa/P_2)$ as a module over $K^T(pt)=\Z[t_1,t_1^{-1},t_2,t_2^{-1}]$. The resulting class is equal to
\begin{multline*}2 \G_{4 1} + 2 \G_{3 2} - \G_{3 3}  - 3 \G_{4 2} + \G_{4 3}\\-\left(t_1+\frac{1}{t_1}+t_2+\frac{1}{t_2}+\frac{t_1}{t_2}+\frac{t_2}{t_1}-6\right)\left(\G_{2 1} - \G_{2 2} - \G_{3 1} + \G_{3 2}\right)\,.\end{multline*}
It can also be written in a shorter form
$$\G_{2 1}\left(2\G_{2 0} -\G_{2 1}-
\left(t_1+\frac{1}{t_1}+t_2+\frac{1}{t_2}+\frac{t_1}{t_2}+\frac{t_2}{t_1}-6\right)(\G_{0} -\G_{1})\right)\,.$$
In terms of Grothendieck roots of $\gamma$ we obtain the solution announced in Theorem \ref{g2alpha}.\qed\medskip

\begin{remark}\rm
For line bundles $u_1$ and $u_2$ we have
$$\G_{ab}(u_1\oplus u_2)=\frac{\left(1-\frac{1}{u_1}\right){}^{a+1}
   \left(1-\frac{1}{u_2}\right){}^b}{1-\frac{u_2}{u_1}}+
   \frac{\left(1-\frac{1}{u_2}\right){}^{a+1}
   \left(1-\frac{1}{u_1}\right){}^b}{1-\frac{u_1}{u_2}}$$
The Grothendieck classes are K-theoretic analogues of the Schur classes in cohomology
$$S_{ab}(u_1\oplus u_2)=\frac{x_1^{a+1}x_2^{b}}{x_1-x_2}+\frac{x_2^{a+1}x_1^{b}}{x_2-x_1}\,,$$
where $x_i=c_1(u_i)$.
The class $\G_{a,b}(u_1\oplus u_2)$ is a polynomial in the inverses of $u_i$'s, which is why we use Grothendieck roots of the dual bundle. For example 
\begin{align*}
 \G_{0}(\gamma^*)= & 1 \\
 \G_{1}(\gamma^*)= & 1-z_1 z_2 \\
 \G_{2}(\gamma^*)= & z_2 z_1^2+z_2^2 z_1-3 z_2 z_1+1 \\
 \G_{11}(\gamma^*)= & \left(1-z_1\right) \left(1-z_2\right) \\
 \G_{3}(\gamma^*)= & -z_2 z_1^3-z_2^2 z_1^2+4 z_2 z_1^2-z_2^3
   z_1+4 z_2^2 z_1-6 z_2 z_1+1 \\
 \G_{21}(\gamma^*)= & \left(1-z_1\right) \left(1-z_2\right)
   \left(1-z_1 z_2\right) \\
\G_{4}(\gamma^*)= & z_2 z_1^4+z_2^2 z_1^3-5 z_2 z_1^3+z_2^3
 z_1^2-5 z_2^2 z_1^2+10 z_2 z_1^2+z_2^4 z_1-5 z_2^3
z_1+\\
&\phantom{xxxxxxxxxxxxxxxxxxxxxxxxxxxxxxxxxxxx}+10 z_2^2 z_1-10 z_2 z_1+1 \\
 \G_{31}(\gamma^*)= & \left(1-z_1\right) \left(1-z_2\right)
   \left(z_2 z_1^2+z_2^2 z_1-3 z_2 z_1+1\right) \\
 \G_{22}(\gamma^*)= & \left(1-z_1\right){}^2 \left(1-z_2\right){}^2
\end{align*}
where $z_1$ and $z_2$ are Grothendieck roots of $\gamma$.
The Chern character transports the Grothendieck classes to Schur classes modulo higher gradation:
$$ch(\G_I(-))=S_I(-)+\sum_{J\supset I}a_JS_J(-)\,.$$
For example of partitions appearing in (\ref{G2nonequiv}):
\begin{align*}ch(\G_{41})&= S_{41}-(S_{42}+2
   S_{51})+\left(\tfrac{1}{2}
   S_{43}+\tfrac{25}{12}
   S_{52}+\tfrac{13}{6}
   S_{61}\right)-\\
&\phantom{xxxxxxxxxxxxxxxxxxxxxxxxx}
   -\left(\tfrac{1}{6}
   S_{44}+\tfrac{13}{12}
   S_{53}+\tfrac{7}{3}
   S_{62}+\tfrac{5}{3}
   S_{71}\right)+\dots,\\
ch(\G_{32})& =S_{32}-
   \left(\tfrac{3}{2}
   S_{33}+\tfrac{3}{2}
   S_{42}\right)+
   \left(\tfrac{7}{3}
   S_{43}+\tfrac{5}{4}
   S_{52}\right)-
   \left(\tfrac{5}{4} S_{44}+2
   S_{53}+\tfrac{3}{4}
   S_{62}\right)+\dots\end{align*}
for bundles of rank two.  The Chern character allows to deduce the formulas in cohomology (\ref{homologicznie}) from those in K-theory. 
\end{remark}

\begin{remark}\rm Formula for the class of $\Gdwa/P_2$ in K-theory is not unique since there are relations between the variables $z_1$, $z_2$, $t_1$ and $t_2$.
When written in terms of Grothendieck classes of $\gamma$ the formula becomes much more complicated\footnote{The nonequivariant class can be written as $-(2 \G_{3 2}(\gamma) + 7 \G_{3 3}(\gamma) + 2 \G_{4 1}(\gamma) + 13 \G_{4 2}(\gamma) +
 37 \G_{43}(\gamma) + 64 \G_{4 4}(\gamma) + 8 \G_{51}(\gamma) + 40 \G_{5 2}(\gamma) +
 106 \G_{53}(\gamma) + 196 \G_{54}(\gamma) + 273 \G_{55}(\gamma))$.}.
The difference appears because, unlike the equality for Schur classes, in general $\G_I(\gamma^*)\not=(-1)^{|I|}\G_I(\gamma)$.
In the case of cohomological fundamental class there is no ambiguity. The classes $S_{50}$, $S_{41}$ and $S_{32}$ span both $\HH^{10}(Gr(2,7))$ and the space of polynomials in two variables of degree 5. Therefore the presentation of $[\Gdwa/P_2]\in \HH^{10}(Gr(2,7))$ in terms of Chern roots of $\gamma^*$ is unique.
\end{remark}

\section{More residue formulas}
As an addendum we list another two residue formulas which are available in literature or just checked by us.
\subsection{Standard  Grassmannian, two sets of characters}
Suppose a bundle $E$ over $Gr(m,n)$, associated to a representation of $GL_m(\C)\times GL_{n-m}(\C)$, is given. Let $f(t_1,t_2,\dots,t_n)$ be its character. It is symmetric with respect to $\Sigma_m\times\Sigma_{n-m}$. Let us write $f$ as a function on two sets of variables $f(A,B)$, $|A|=m$, $|B|=n-m$. It is symmetric with respect to $A$-variables and $B$-variables.  Let
\begin{itemize}
\item $T=\{t_1,t_2,\dots, t_n\}$

\item $Z_1=\{z_1,z_2,\dots, z_m\}$

\item $Z_2=\{z_{m+1},z_{m+2},\dots ,z_n\}$

\item $Z=Z_1\cup Z_2$
\end{itemize}
Then
\begin{align*}p_!(E)&=\sum_{A\subset T,\;|A|=m}\frac{ f(A,\oA{^{-1}})}{[\oA/A]}\\
&=\frac{1}{m!(n-m)!} Res_{Z=0,\infty} \frac{f(Z_1,Z_2)\cdot[R(Z_1)]
\cdot[Z_1/Z_2]\cdot [R(Z_2)]
\cdot \dlog(Z)}{[T/Z]}\,.\end{align*}
Another formula is given in  \cite[Prop 5.1]{All}:
$$p_!(E)= Res_{Z=0,\infty} \frac{f(Z)\cdot[R^+(Z)]\cdot \dlog(Z)}{[T/Z]}\,.$$
Here there is no factor $\frac1{m!(n-m)!}$, and the numerator is much simpler.

\subsection{Flag manifold}
In \cite[Prop 5.3]{All} Allman gives the following formula for the push-forward on the full flag manifold.
Let $Fl(n)=GL_n(\C)/B$ be the full flag manifold. Let $E$ be a vector bundle over $Fl(n)$ associated to a representation of the maximal torus. Let $f\in K^\T(pt)$ be the class of this representation.
Let
\begin{itemize}
\item $T=\{t_1,t_2,\dots, t_n\}$

\item $Z=\{z_1,z_2,\dots, z_n\}$
\end{itemize}
Then the formula for the push-forward is the following:
$$p_!(E)=\sum_{\sigma\in\Sigma_n}\frac{ f(\sigma T)}{[R^+(\sigma T^{-1})]}=
 Res_{Z=0,\infty} \frac{f(Z)\cdot [R^+(Z)]\cdot\dlog(Z)}{[T/Z]}\,.$$
In particular if $$f(T)=
\prod_{i=1}^n\left(1-\frac1{t_i}\right)^{a_i+n-i}$$ we obtain the formula from \cite[Prop. 3.3]{RS} for Grothendieck polynomial $$G_{(a_1,a_2,\dots,a_n)}(t_1,t_2,\dots, t_n)\,.$$

\subsection{$\Gdwa/B$}
The formula (\ref{g2alpha}) is valid for $\Gdwa/B$ because the  $\Gdwa/B$ embeds into $Fl_{1,2}(\C^7)$ and the class of $\Gdwa/B$ comes from the Grassmannian $Gr(2,7)$ via  pull-back:
$$[\Gdwa/B]=\pi^*[G/P_2]\,,$$
where $\pi\colon Fl_{1,2}(\C^7)\to Gr(2,7)$ is the natural projection.

\subsection{Quadric}\label{kwadryka}
The homogeneous space for $O(2n)$ consisting of isotropic lines $IG(1,2n)\subset \PP^{2n-1}$ does not fit to the previous schemes. Suppose $E$ is a vector bundle over $IG(1,2n)$ which is associated to a representation of $\C^*\times O(2n-2)$. Its character $f\in K^\T(pt)$ is symmetric with respect to permutations of $t_2,t_3,\dots,t_n$ and with exchanging $t_i$ with $t_i^{-1}$ for $i\geq 2$.
Let
\begin{itemize}
\item $T=\{t_1,t_2,\dots, t_n\}$

\item $Z_2=\{z_{2},z_{m+2},\dots ,z_n\}$

\item $Z=\{z_1\}\cup Z_2$

\end{itemize}

One has
\begin{multline*}p_!(E)=\sum_{a\in T^{\pm  1}}\frac{f(a,T\setminus \{a\})}
{[(T^{\pm 1}\setminus \{a,a^{-1}\})/\{a\}]}=\\ \frac1{2^{n-1}}Res_{Z=0,\infty}
\frac{f(z_1,Z_2)\cdot[\{z_1 ^{ - 2}\}]\cdot  [Z ^{ - 1}\otimes Z_2 ^{ - 1}]\cdot [R^+(Z)]\cdot  \dlog(Z)}{[T^{\pm1}/ Z]}\end{multline*}

This case is really different from the previous ones. It cannot be reduced to the push-forward of classes expressed by symmetric function in groups of variables. Here the second set of variables is invariant with respect to the action of $\Z_2^{n-1}\rtimes \Sigma_{n-1}$. We leave this formula as an exercise for the reader. Our proof is just by an elementary calculus of functions.

\section*{Appendix: About computations}

Let us briefly discuss the methods of computations.
To find the push-forward of a class in K-theory given by the polynomial $f(t_1,t_2)$ we apply the localization formula (\ref{locformula}). The Weyl group of $\Gdwa$ is isomorphic to the dihedral group $D_{12}$. The fixed points of $\Gdwa/P_2$ are indexed by the cosets $D_{12}/\Z_2$.  The contribution coming from the coset of the identity is equal to
$$\theta(t_1,t_2)=\frac{f(t_1, t_2)}{\big(1 - t_1\big)  \big(1 - t_2\big) \big(1 - t_1 t_2\big) \big(1 - \frac{t_1^2}{t_2}\big)\big(1 - \frac{t_2^2}{t_1}\big)}\,.$$
The remaining fixed points lie in the orbit of $\Z_6$. The action of the generator of $\Z_6$ is given by the substitution:
$$\xi:(t_1,\,t_2)\;\mapsto \;(t_2,\, t_2/t_1)\,.$$
The push-forward $p_{\Gdwa!}$ is equal to the sum of six quotients:
\begin{equation}\label{cyclicsum}\sum_{k=0}^5\theta(\xi^k(t_1,t_2))\,.\end{equation}
We find that
the push-forwards of the Grothendieck polynomials for the partitions
$$
 [0],\; [1],\; [2],\; [11],\; [3],\; [21],\; [31],\; [22]
$$
are equal to one. Moreover
$$p_{\Gdwa!}\left(\G_{4}\right)=p_{\Gdwa!}\left(\G_{41}\right)=
 p_{\Gdwa!}\left(\G_{32}\right)=2\,,\qquad p_{\Gdwa!}\left(\G_{33}\right)=0\,.$$
The remaining values of $p_{\Gdwa!}\left(G_{I}\right)\in K^\T(pt)$ can be expressed in variables
$$A=3-t_1-\frac{1}{t_2}-\frac{t_2}{t_1}\,,\qquad B=3-t_2-\frac{1}{t_1}-\frac{t_1}{t_2}$$
\def\sskip{\noalign{\vskip 2mm}}
$$\begin{array}{l}
 p_{\Gdwa!}\left(\G_{5}\right)=
 p_{\Gdwa!}\left(\G_{51}\right)=
 p_{\Gdwa!}\left(\G_{42}\right)=A+B\, \\ \sskip
 p_{\Gdwa!}\left(\G_{52}\right)=A^2+B^2+A B-2 A-2 B\, \\ \sskip
 p_{\Gdwa!}\left(\G_{43}\right)=A B-A-B\, \\ \sskip
 p_{\Gdwa!}\left(\G_{53}\right)=A B (A+B-5)\, \\ \sskip
 p_{\Gdwa!}\left(\G_{44}\right)=-A B\, \\ \sskip
 p_{\Gdwa!}\left(\G_{54}\right)= A^2 B^2+A^3+B^3-3 A^2 B-3 A B^2-3 A^2-3 B^2+8 A B \,\\ \sskip
 p_{\Gdwa!}\left(\G_{55}\right)=-A^2 B^2-A^3-B^3+2 A^2 B+2 A B^2+2 A^2+2 B^2-4 A B\, \\
\end{array}
$$
Next we compute the intersection matrix
$$m_{I,J}=p_!(\G_I\cdot \G_J)$$
for the partitions $I$ and $J$ contained in the $5\times 2$ rectangle. This is done by the application of the formula (\ref{locformulagrass}). We obtain a matrix of the size $21\times 21$. Its entries are Laurent polynomials in $t_1$ and $t_2$. The entries again can be expressed in $A$ and $B$ as symmetric polynomials. The matrix after reduction of variables is still quite large, hence we do not give its precise form, see \cite{Wolfram}. It is remarkable that the determinant of this matrix is equal to $-1$.
We invert the intersection matrix and we find that
\begin{multline*}p_{\Gdwa!}(\G_I)=
p_!\Big(
\G_I\cdot \big( 2 \G_{4 1} + 2 \G_{3 2} - \G_{3 3}  - 3 \G_{4 2} + \G_{4 3}+\\+(A+B)\left(\G_{2 1} - \G_{2 2} - \G_{3 1} + \G_{3 2}\right)\big)\Big)\,.\end{multline*}
for all partitions $I$. It follows that
\begin{multline*}[\Gdwa/P_2]=  2 \G_{4 1} + 2 \G_{3 2} - \G_{3 3}  - 3 \G_{4 2} + \G_{4 3}+\left(A+B\right)\left(\G_{2 1} - \G_{2 2} - \G_{3 1} + \G_{3 2}\right)\,.\end{multline*}

\end{document}